\documentclass[12pt]{article}

\usepackage{thmtools}
\usepackage{thm-restate}
\usepackage[pdftex,pdfstartview=FitH,pdfborderstyle={/S/B/W 1}, colorlinks=true,linkcolor=blue,urlcolor=blue,citecolor=blue]{hyperref}
\usepackage{cleveref}

\usepackage{filecontents}
\usepackage{ifpdf}
\usepackage{latexsym}
\usepackage{amssymb}
\usepackage{pifont}
\usepackage{amsopn}
\usepackage{graphicx}
\usepackage{tikz}

\usepackage[empty]{fullpage}

\usepackage{bookmark}

\newtheorem{theorem}{Theorem}[section]
\newtheorem{proposition}[theorem]{Proposition}
\newtheorem{lemma}[theorem]{Lemma}
\newtheorem{corollary}[theorem]{Corollary}
\newtheorem{definition}[theorem]{Definition}
\newtheorem{notation}[theorem]{Notation}

\newcommand{\zz}{\mathbb{Z}}
\newcommand{\rr}{\mathbb{R}}
\newcommand{\cc}{\mathbb{C}}
\newcommand{\tH}{\mathrm{th}}
\newcommand{\nn}{\mathbb{N}}

\newcommand{\diam}{\mathrm{diam}}

\begin{document}

\pagestyle{plain}

\title{Meromorphic $L^2$ functions on flat surfaces}
\author{Ian Frankel}

\maketitle

\begin{abstract}We estimate spectral gaps for the Hodge norm on quadratic differentials. To each tangent direction at any point $(X,q)$ in the principal stratum of quadratic differentials, we associate a Hodge norm, and control the logarithmic derivative of vectors perpendicular to the principal directions in terms of the $q$-areas of the components corresponding to thick-thin decompositions and the lengths of short curves in the $q$-metric. In the worst case scenario, one gets a spectral gap of size $C_{g,n}\mathrm{sys}(X,q)^2$. \end{abstract}

\tableofcontents

\section{Introduction}

\noindent In \cite{Hypertrain}, we defined a metric $d_E$ on the space of unit area quadratic differentials on Riemann surfaces in terms of flat geometrty, and wrote $d^{ss}$ for the restriction of this metric to strongly stable leaves of the Teichm\"uller geodesic flow. We showed this flow exhibits hyperbolic behavior near quantitatively recurrent trajectories. More precisely,

\begin{theorem}\label{DefiniteContraction}Let $K$ be a compact subset of the unit tangent space to the moduli space $\mathcal{M}_{g,n}$, and let $\lambda$ denote Lebesgue measure on $\rr$. Let $\theta \in (0,1)$. There are positive real constants $C(K), r_0(K), \alpha(K,\theta)$ such that the following holds: if $T > 0,$ and $X, g_T(X) \in K$ and moreover, $$\lambda \{t \in [0,T]: g_t(X) \in K\} > \theta T, ~ \mathrm{and}$$ $$d^{ss}(X_1,X),d^{ss}(X_2,X) < r_0(K), ~ \mathrm{then}$$ $$\frac{d^{ss}(g_T X_1, g_T X_2)}{d^{ss}(X_1,X_2)} < C(K)e^{-\alpha(K,\theta)T}.$$ \end{theorem}

\noindent A weaker version of this, using a construction called the \emph{modified Hodge norm}, was used in \cite{abem}, \cite{EM}, and \cite{HB} to provide Teichm\"uller space analogues to the results in the thesis of Margulis \cite{Margulis} on the geodesic flow of compact negatively curved manifolds.\\

\noindent A substantial difficulty in the implementation of this technique came from the fact that the Hodge norm arguemnts depended on a compactness argument, and estimates from this argument became degenerate as quadratic differentials approached the boundary of the principal stratum. The purpose of \Cref{DefiniteContraction} was to remove this technical difficulty by giving an alternative compactness argument that did not degenerate at the multiple zero locus.\\

\noindent The purpose of this paper is to give a more accurate description of the hyperbolic behavior of the Hodge norm that does not depend solely on compactness arguments. In addition, the estimates persist as one approaches the boundary of a stratum, and we estimate the rate at which hyperbolicity is lost as one approaches the boundary of moduli space; in some instances we even see long stretches in the cusps with no loss of hyperbolicity. In particular, the Hodge norm exhibits hyperbolic behavior over stretches when the surface does not make progress in more than one component of a Minsky product region. (See \cite{ProdReg} or \cite{Fellow} for precise descriptions of product regions.)\\

\noindent In particular, we define a number called the \emph{Hourglass Ratio} $H(X,q)$ of a half-translation surface $(X,q)$, which is equal to 1 if $X$ is in the thick part of $\mathcal{M}_{g,n}$, and in the cusps is less than or equal to $\ell/\sqrt{A}$ if the surface can be cut into two essential or cylinder pieces with area at least $A$ by a system of curves of length $\ell$ in the $q$-metric.\\

\noindent Our main theorem is the following:\\

\begin{theorem}\label{HodgeGap}Let $(X,q)$ be any quadratic differential in the principal stratum of half translation surfaces of genus $g$ with $n$ marked points. Let $(\tilde{X},\omega)$ be the orienting double cover for $(X,q)$. Let $g_t$ denote the Teichm\"uller flow, with $g_t(X,\omega) = (X_t,\omega_t)$. Let $\|\alpha\|_t$ denote the Hodge norm of $\alpha \in H_{odd}^{1}(\tilde{X}),$ where $\alpha$ is the cohomology class of a (-1)-eigenform of the involution of $X$ respecting the double cover, which is orthogonal to $\mathrm{Re}(\omega)$ and $\mathrm{Im}(\omega)$ with respect to the Hodge norm on $H^1(\tilde{X_t},\rr)$. Then $$\frac{d}{dt}\left|\log(\|\alpha\|_t)\right| < 1 - C_{g,n}H(X,q)^2.$$ \noindent where $C_{g,n} > 0$ is a constant depending on only $g$ and $n$.\end{theorem}

\noindent Remark: This provides an infinitesimal spectral gap result for a norm on the tangent space at almost every point. The estimate does not degenerate near the multiple zero locus, even though the norm does. Thus we can give a qualitative estimate of how the gap decays in the cusp of the moduli space.\\

\noindent To make the analogy to negatively curved Riemannian manifolds, parallel transport of geodesic flow contracts perpendicular tangent directions at an exponential rate, where the exponent is the proportional to the square root of the sectional curvature. Although the Teichm\"uller metric is not Riemannian, this gives an estimate akin to ``the sectional curvature is at most $-C_{g,n}H(X,q)^4$ in all 2-planes containing the direction of $q$," at least for the purposes of estimate expansion/contraction of the unstable/stable manifolds for the flow, with the caveats that must be applied whenever attempting to use the Hodge norm as a norm on the tangent space to $\mathcal{M}_{g,n}$.\\

\noindent As an immediate corollary, we have

\begin{theorem} Let all notation be as in theorem 1. Let $(X,q)$ be a unit area quadratic differential in the principal stratum, whose shortest essential simple closed curve has length $Sys(X,q)$ in the $q$-metric. Then $$\frac{d}{dt}\left|\log(\|\alpha\|_t)\right| < 1 - C_{g,n}Sys(X,q)^2.$$
\end{theorem}

\noindent Proof: This is immediate since the systole is less than or equal to the hourglass ratio (up to a constant multiple). $\Box$

\subsection{Acknowledgments} The author would like to thank Alex Eskin, Jeremy Kahn, Kasra Rafi, and Alex Wright for conversations that provided the inspiration for this project. The author was supported by a Fields Postdoctoral Fellowship and a research fellowship from National Research University Higher School of Economics.

\includegraphics{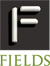} \includegraphics{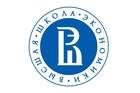}

\section{Quadratic Differentials and the Hodge Norm}

\subsection{Quadratic Differentials}

\noindent Let $X$ be a Riemann surface of genus $g$ with $n$ punctures. Assume $X$ admits a conformal metric of constant curvature $-1$ with finite area, i.e. $3g-3+n>0$. A \emph{quadratic differential} $q$ on $X$ is a holomorphic section of the tensor square of the cotangent bundle of $X$ that extends meromorphically to the compact Riemann surface $\overline{X}$ obtained by adding in a point at each cusp, such that there are no poles of order 2 or higher. The space of quadratic differentials on $X$ has $3g-3+n$ complex dimensions.\\

\noindent Such $(X,q)$ admits a system of holomorphic charts covering all points of $X$ at which $q$, such that if $z$ is the local coordinate then $q = dz^2 := dz \otimes dz$. Call these coordinate charts $q$-\emph{charts}. The change of charts preserve a metric as well as foliations by vertical and horizontal lines. The metric and foliations have singular extensions to all of $\overline{X}$ with cone-type singularities angle $(2+k)\pi$, and $(2+k)$-pronged singularities, at each point where $q$ vanishes to order $k$, where $k \geq -1$. We will mostly concern ourselves with the generic case, in which $q$ has $n$ simple poles and $4g-4+n$ simple zeros. Such $(X,q)$ is said to belong to the \emph{principal stratum}. We will also deal with quadratic differentials that are squares of holomorphic $1$-forms.

\subsection{Teichm\"uller Deformations}

\noindent For $t \in \rr$ let $g_t$ be the matrix $\left(\begin{array}{cc} e^t & 0\\ 0 & e^{-t}\end{array}\right),$ acting on $\cc = \rr^2$ with the usual identification $\left(\begin{array}{c} x\\ y\end{array}\right) = x + iy$. Then $(X,q)$ is given by a system of charts with transition maps of derivative $\pm 1$, and $g_t$ simultaneously acts on these charts. Let $(X_t,q_t)$ denote the 1-parameter family of Riemann surfaces of genus $g$ with $n$ marked points whose complex structures are given by images of the $q$-charts under $g_t$ in $\rr^2$. (This respects the identifications). It is well known that $X_t$ is a geodesic parametrized by arc length with respect to the Teichm\"uller metric on the moduli space of Riemann surfaces.

\subsection{Orienting Double Cover}

\noindent Let $(X,q)$ be in the principal stratum. The \emph{orienting double cover} of $(X,q)$ is the Riemann surface obtained by taking the branched cover of $\overline{X}$ which is branched of order $2$ over all zeros and poles of $q$, and such that the preimages of the nonsingular points of the $q$-metric are given by vertical unit tangent vectors with respect to the $q$-metric.\\

\noindent On $\tilde{X}$, the pullbacks of the vertical and horizontal foliations of $q$ are orientable, and induced by the square of a holomorphic $1$-form $\omega$ on $X$. The induced double cover therefore admits a collection of charts whose transition maps are translations, i.e. the structure of a Riemann surface with a holomorphic $1$-form. This $1$-form is a $-1$-eigenvector for the holomorphic involution $\iota: \tilde{X} \to \tilde{X}$ that permutes the sheets of the double cover. The $1$-form $\omega$ has $4g-4+n$ double zeros, which are all fixed by $\iota.$ As is well known, the locus of holomorphic $1$-forms realized by these double covers is locally described by the $-1$-eigenspace the action $\iota$ on the first de Rham cohomology of $\tilde{X}$, which we will call $H_{odd}^1(\tilde{X},\cc)$.\\

\noindent Since $\iota$ is holomorphic it commutes with the Hodge star operator, and the real and imaginary parts of a holomorphic 1-form $\alpha + i \beta$, $\alpha,\beta$ real, are related by Hodge star: $\beta = *\alpha$, $\alpha = -*\beta$. Thus $H_{odd}^1(\tilde{X},\rr)$ consists of those classes representable as real parts of holomorphic $1$-forms in $H_{odd}^1(X;\cc)$. Let $\eta$ be an odd holomorphic $1$-form. Then $(\iota/\eta)^2$ is $\iota$-invariant and therefore descends to a meromorphic function $f$ on $X$. A dimension count in fact shows that $f$ is of the form $q_2/q$, where $q_2$ is another quadratic differential on $X$. (For any $q_2$ one gets such a function.)

\subsection{Hodge Norm and its First Variation}

\noindent The \emph{Hodge Norm} of $c \in H^1(M,\rr)$ on a compact Riemann surface $M$ is given by $\|c\|_X^2 = \int_M \mathrm{Re} (h_c) \wedge \mathrm{Im} (h_c)$ where $h_c$ is the holomorphic $1$-form whose real part is in the class $c$. The Hodge norm comes from a Hilbert space inner product:

$$\langle c_1, c_2\rangle_M = -\frac{1}{2}\mathrm{Im}\left(\int_X h_{c_1} \wedge \bar{h}_{c_2}\right).$$

\noindent As the complex structure of $M$ varies continuously, one can think of this group as a fixed vector space with a varying norm.\\

\noindent Let $(M_t,q_t)$ is a Teichm\"uller geodesic given by a 1-parameter family from a quadratic differential $q_0 = \omega_0^2$ defined on $M = M_0$, where $\omega_0$ us a holomorphic 1-form, and assume that the $q$-metric has unit area. There is a unit area holomorphic 1-form $\omega_t$ on $M_t$, whose cohomology class is $e^t \mathrm{Re}([\omega_0]) + i \cdot e^{-t}\mathrm{Im}([\omega_0])$. The area form $|q|$ is given by $\mathrm{Re}(\omega) \wedge \mathrm{Im}(\omega).$ If $c \in H^1(M,\rr)$, Forni's variational formula (\cite{FE} lemma $2.1^\prime$, see also \cite{FMZ}, section 2.6) says that $$\left.\frac{d}{dt} \left(\|c\|_{M_t}^2\right) \right|_{t = 0} = -2\mathrm{Re}\left(\int_{M} \left(\frac{h_c}{\omega}\right)^2 |q|\right).$$

\noindent One sees that in the simplest case, $h_c = \pm \omega$, the derivative of the Hodge norm squared is -2, since one may apply the matrices $\left(\begin{array}{cc} 1 & 0\\ 0 & e^{-2t}\end{array}\right)$ to obtain the family $M_t$, and the translation surface structures on $M_t$ whose holomorphic $1$-forms have cohomologous real parts. It follows that the logarithmic derivative of the Hodge norm of $c$ is $-1$. Similarly, if $h_c = \pm i\omega$ one gets that the logarithmic derivative of the Hodge norm is $1$.\\

\noindent The variational formula tells us that these are the two extreme behaviors. A simple application of the Cauchy-Schwarz inequality implies that the logarithmic derivative of the Hodge norm is always in $[-1,1],$ and the lower and upper extreme cases are attained by real and imaginary multiples of $\omega$ respectively. To attain statements about uniform hyperbolicity on certain regions of the moduli space, we would like to attain strict bounds for the growth rate of cohomology classes on the subspace $$\langle[\mathrm{Re}(\omega_t)],[\mathrm{Im}(\omega_t)]\rangle^\perp \cap H_{odd}^1(X;\rr),$$ i.e. to bound the logarithmic growth rate away from $\pm 1$. Here, the symbol $\perp$ is the symplectic complement with respect to the cup product. As it happens, it is also the Hodge orthogonal complement. This subspace is $g_t$-invariant.\\

\noindent The way we will proceed is as follows: Suppose that $\tilde{f}$ is the pullback of $f = q_2/q$ to the orienting double cover of $(X,q)$, where $q,q_2$ are quadratic differentials on $X$ and $\|\tilde{f}\|_2 = 1$. Let $(M_t,\omega_t) = g_t(\tilde{X},\omega)$ be the orienting double covers corresponding to the family $g_t(X,q)$. Write $|\omega^2|$ for the area form obtained by pulling back $|q|$. Let $c \in H^1(\tilde{X};\rr)$ be such $h_c/\omega$ is the pullback of $q_2/q$. Then if $\|\mathrm{Im}(\tilde{f})\|_2 = \delta$, we have $$\begin{array}{rl}\left.\frac{d}{dt} \left(\|c\|_{M_t}^2\right) \right|_{t = 0} & = -2\mathrm{Re}\int_{M} (\tilde{f})^2 |\omega^2|\\ {} & = -2\int_M \mathrm{Re}(\tilde{f})^2 - \mathrm{Im}(\tilde{f})^2|\omega^2| \\ {} & = -2(1 - \delta^2 - \delta^2).\end{array}$$

\noindent It follows that $$\left.-\frac{d}{dt}\log(\|c\|_{M_t})\right|_{t = 0} = 1 - 2\delta^2.$$

\noindent We have a similar analysis when $\|\mathrm{Re}(\tilde{f})\| = \delta$. In particular, if $\|\mathrm{Re}(\tilde{f})\|, \|\mathrm{Im}(\tilde{f})\| \geq \delta,$

$$\left|\frac{d}{dt}\log(\|c\|_{M_t})\right|_{t = 0} \leq 1 - 2\delta^2.$$

\noindent In the interest of proving spectral gap results, it is enough to prove $$\|\mathrm{Im}(\tilde{f})\|_2/\|\tilde{f}\|_2 = \|\mathrm{Im}(f)\|_2/\|f\|_2 \geq \delta$$ whenever $(X,q)$ satisfies some geometric conditions on $(X,q)$ depending on $\delta$. The equality is automatic. From now on, we will assume $\|f\|_2 = 1$ and that $f$ is meromorphic and orthogonal to the constants in $L^2(X,|q|).$ These will be all we need to assume about $f$, and in the case when $q$ is in the principal stratum, such $f$ consist of exactly the functions $q_2/q$ with $\int_X (q_2/q)|q| = 0,$ up to scale.

\section{Preliminary notation and estimates}

\subsection{Some Fixed Notation}

\noindent Throughout, $C,C^\prime,C^{\prime\prime}$ will be topological constants, whose values may vary from one line to the next.\\

\noindent Throughout, we will use the symbol $u ~ \dot{\asymp} ~ v$ to mean that there is some constant $C \in (1,\infty)$, depending only on the particular moduli space we are considering, such that $$C^{-1} < u/v < C,$$ where $u,v$ are two positive quantities that depend on a point in the moduli space.\\

\noindent $\mathbf{1}$ denote the constant function everywhere equal to $1$; let $\mathbf{1}_V$ denote the characteristic function of a set $V$.

\subsection{A Gradient Estimate}

\noindent Let $(X,q)$ be a point in quadratic differential space. For a measurable function $f: X \to \cc$, let $\|f\|_p$ denote the $L^p$ norm of $f$ with respect to the area measure $|q|$ of the $q$-metric. Let $\langle f,g \rangle = \int_X f \bar{g} |q|$.\\

\noindent We have the following estimate by Trevi\~no for the gradient of $f$ at $x$ with respect to the $q$-metric. What we state here is a trivial generalization of equation (26) in \cite{Rodrigo}, which requires only that $f$ is holomorphic on a round disk centered at $x$.

\begin{proposition}\label{Rodrigo} Let $\Sigma$ denote the set of singular points of $(X,q)$. Let $B_R(x)$ denote the $R$-neighborhood of $x \in X$ with respect to the $q$-metric. Let $r(x)$ be the largest $R$ such that $B_R(x)$ is an embedded disk not meeting $\Sigma$. For any fixed $t \in (0,1]$, there is a constant $c_t$ such that $$|\nabla f(x)| \leq c_t \|\mathrm{Im}(f) \cdot \mathbf{1}_{B_{t \cdot r(x)}(x)}\|_2.$$\end{proposition}

\subsection{Holomorphic functions on expanding annuli}

\noindent The following will be called the \emph{Small Norm on Expanding Annuli Lemma}. It applies to expanding annuli, which we will define below:

\begin{definition} An expanding annulus $A$ is an open annulus in a surface with a flat metric, following the following additional conditions:

\begin{itemize} \item The metric completion $\bar{A}$ of $A$ as a path metric space has two boundary components $\gamma_1$ and $\gamma_2$, which are rectifiable curves, at most one of which is degenerate (length 0).

\item There is a number $W$, called the \emph{width} of $A$ such that every point on $p \in \bar{A} \setminus A$ is distance $W$ from the boundary component not containing $p$.

\item $A$ is foliated by \emph{level circles} which are differentiable, piecewise smooth, and \emph{monotonically curved}, and each level circle has nonzero geodesic curvature at some point.\end{itemize}

\noindent Here, a \emph{level circle} is a topological circle in $A$, all of whose points are some distance $d$ and $W-d$ from the two boundary components, with $0 < d < W$. \emph{Monotonically curved} means that for some orientation, the geodesic curvature is always non-negative.\end{definition}

\noindent The condition that geodesic curvature is not everywhere zero means that the circumference of a level circle is strictly monotonic in $d$; the annulus is said to be \emph{expanding} in the direction in which the circumferences increase. For the purpose of this lemma we will also want to assume that our expanding annuli are free of singularities. We will choose specific expanding annuli later to tackle the general case, but we can already state the lemma in enough generality here.\\

\begin{lemma}\label{shells} Let $W$ be the union of regions $W_j, j_0 \leq j \leq N$, and assume there are additional regions $W_{{j_0} - 1}$ and $W_{N+1}$, such that the following geometric conditions hold:

\begin{itemize}
    \item $\leq C \cdot 2^{-j} \leq \diam_q(W_j) \leq C^\prime 2^{-j}$
    \item $\leq C \cdot 2^{-2j} \leq Area(W_j) \leq C^\prime 2^{-2j}$
    \item The $2^{-j-2}$-neighborhood of $W_j$ is contained in $W_{j-1} \cup W_j \cup W_{j+1}$ for $j_0 \leq j \leq N$
\end{itemize}

\noindent Suppose, that $f$ is holomorphic on $W$, $\|f\|_2 \leq 1$, and moreover, that $\mathrm(i)~\|\mathbf{1}_{W_{j_0}} \cdot f\|_\infty < \delta_0 2^{j_0}$, and $\mathrm(ii)~\|\mathrm{Im}(f)\|_2 < \delta_0.$ Then $$\mathrm{(a)}~\|\mathbf{1}_{W_j}f(x)\|_\infty < C \delta_0 2^{j},~\mathrm{(b)}\|\mathbf{1}_W\cdot f\|_1 < C\delta_0 2^{j_0},~\mathrm{and}~\mathrm{(c)}~\|\mathbf{1}_W \cdot f\|_2 < C\delta_0.$$ \end{lemma}

\noindent Proof: First, write $I_j$ for $\|\chi_{W_j} \cdot \mathrm{Im}(f)\|_2$. By \Cref{Rodrigo} estimates, and the fact that the area of $W_j$ is comparable to $2^{-2j}$ we have $$\|\mathbf{1}_{W_j}|\nabla f(x)|\|_\infty \leq  C 2^{2j}(I_{j-1}+I_j+I_{j+1}).$$

\noindent Also write $M_j = \|\mathbf{1}_{W_j} \cdot f\|_\infty.$ We have $M_j - M_{j+1} \leq C 2^{j}(I_{j-1}+I_j+I_{j+1})$. Let $j_0$ be the smallest $j$ for whichever annulus we are estimating. Since $\sum\limits_j I_j^2 \leq \delta_0^2,$ the Cauchy-Schwarz inequality gives (a) by $$|M_j - M_{j_0}| \leq \sum_{k = j_0}^j C2^k(I_{k-1}+I_j+I_{k+1}) \leq C^\prime 2^j \delta_0.$$

\noindent For convenience we will write $I_{k_*}$ for $(I_{k-1}+I_k+I_{k+1})$; we have $\sum_k I_{k_*}^2 < 9\delta_0^2.$\\

\noindent There might be infinitely many $W_j$ in $W$, but we may assume there are finitely many and apply the monotone convergence theorem. That is, we show that $\|\mathbf{1}_{W_N^*}\cdot f\|_2 < C\delta_0^2,$ where $W_N^*$ is any finite union $W^* = \bigcup\limits_{j = j_0}^N W_j$, such long as the constant $C$ in the bound does not depend on $N$.\\

\noindent Part (b) comes from summing part (a) times the area of $W_j$ over all $j$.\\

\noindent Proof of (c):

$$\begin{array}{rl}\|\mathbf{1}_W \cdot f\|_2^2 & \leq \sum\limits_{j = j_0}^{N} Area(W_j)M_j^2\\ {} & \leq \sum\limits_{j = j_0}^{N} C2^{-2j}M_j^2\\ {} & \leq C\sum\limits_{j = j_0}^{N} 2^{-2j} \left[M_{j_0}+\sum\limits_{k = j_0}^j C^\prime 2^{k} I_{k_*}\right]^2\\ {} & \leq C\sum\limits_{j = j_0}^{N} 2^{-2j}M_{j_0}^2 + \left[\sum\limits_{k = j_0}^j C^\prime 2^{-j}2^k I_{k_*}\right]^2\\ {} & \leq C(2^{-2j_0}M_{j_0}^2) + C^\prime \sum\limits_{j = j_0}^{N}\left[\sum\limits_{k = j_0}^j C^\prime 2^{k-j} I_{k_*}\right]^2\\ {} & \leq C\delta_0^2 + C\sum\limits_{j=j_0}^N \sum\limits_{k=j_0}^j\sum\limits_{r=j_0}^j 2^{k-j}2^{r-j} I_{k_*}I_{r_*}\\ {} & \leq C\delta_0^2 + C\sum\limits_{j=j_0}^N \sum\limits_{k=j_0}^j 2^{k-j}\sum\limits_{r=j_0}^j 2^{r-j} (I_{k_*}^2 + I_{r_*}^2)\\ {} & = C\delta_0^2 + C\sum\limits_{j=j_0}^N \sum\limits_{k=j_0}^j 2^{k-j}\sum\limits_{r=j_0}^j 2^{r-j} 2I_{k_*}^2\\ {} & \leq C\delta_0^2 + C\sum\limits_{j = j_0}^N \sum\limits_{k = j_0}^j 2^{k-j}2I_{k_*}^2 \\ {} & \leq C \delta_0^2 + C \sum\limits_{k = j_0}^N I_{k_*}^2 \sum\limits_{j = k}^N 2^{k-j}\\ & \leq C \delta_0^2 + C \sum\limits_{k = j_0}^N I_{k_*}^2 \\ {} & \leq C \delta_0^2.\end{array}$$

\noindent Now take square roots. $\Box$\\

\noindent We would like to point out a simple intuitive reason to expect a lemma like this: consider any Laurent series convergent on the region $\epsilon < z \leq 1/\epsilon$ extending continuously to the boundary with small values of $z$ on the outer boundary. The constant term positive powers of $z$ must have small coefficients, and the negative powers of $z$ all contribute roughly the same to the real and imaginary parts on any circle $|z| = R.$

\subsection{Holomorphic functions on cylinders}

\noindent In this section we examine the case in which our surface has a cylinder of large modulus. We define cylinder and modulus:

\begin{definition} Let $(X,q)$ be a Riemann surface with quadratic differential. A \emph{cylinder} in the $q$-metric is an open annulus $A$ in $X$ which is the disjoint union of all singularity-free $q$-geodesic representatives of a non-trivial free isotopy class $\alpha$ of simple closed curves in $X$. The curve $\alpha$ is said to be a \emph{cylinder curve} for $(X,q)$.\end{definition}

\begin{definition}The \emph{modulus} of a cylinder with core curve $\alpha$ on a surface $X$ with quadratic differential $q$ is given by $\frac{w_q(\alpha)}{\ell_q(\alpha)}$, where $\ell_q(\alpha)$ is the length of a $q$-geodesic representative of $\alpha$, and $w_q(\alpha)$ is the length of the shortest arc that passes through every $q$-geodesic representative of $\alpha$.\end{definition}

\noindent If a Riemann surface $Y$ is homeomorphic to an annulus, it admits a singularity-free quadratic differential metric for which the whole of $Y$ is a cylinder, and this quadratic differential is unique up to scalars, so the metric is unique up to homothety. The \emph{modulus of the annulus} of the annulus is modulus of this cylinder. If the modulus is infinite, then there are two conformal classes: quotient of the upper half-plane by $(\zz,+)$, which we will call \emph{semi-infinite} and the other by the quotient of $\cc$ by $(\zz,+)$, which we call \emph{bi-infinite}. The only connected Riemann surfaces properly containing bi-infinite annuli are $\cc$ and $\hat{\cc}$ so we will never encounter these. See e.g. \cite{Hubbard} for detail on the uniformization of annuli.\\

\noindent Fix $b < 1$. Write $B$ for $b^{-1}$. Pick a local coordinate system on a cylinder $(-L,L) \times [0,2\pi s]$ with the identification $(x,0) \thicksim (x,2\pi s)$ for all $x \in -[L,L].$ We will assume $L > 2b$.\\

\noindent The restriction of $f$ to $A^*$ can be viewed as a function on the annulus $|\log(|z|)| < L$ via the map $(x,y) \mapsto z = e^{(x+iy)/s}.$ Since $f$ is holomorphic on this annulus, it is given by a convergent Laurent series $f(z) = \sum\limits_{n = -\infty}^\infty \hat{f}(n)z^n$. For any $\epsilon \in (0,L/s),$ the sum converges absolutely and uniformly on the annulus $|\log(|z|)| < L/s - \epsilon$.\\

\noindent We would now like to make the following important estimate, which is true in general for any cylinder with coordinates $(-L,L) \times [0,2\pi s]/(x,0) \thicksim (x,2\pi s),$ and any $f$ which is holomorphic on the cylinder and whose imaginary part has small $L^2$ norm.

\begin{lemma}\label{CylinderConstant} Fix a small number $b \in (0,1/2)$, and set $B = b^{-1}$. $f(x,y) - \hat{f}(0).$ Assume $f$ is meromorphic on a cylinder whose core curve has length $2\pi s$ and whose width is $2L$. There is a constant $C$, depending only on $b$, such that the following hold whenever $L/s$ is sufficiently large:\\

\noindent $\mathrm{(a)}~|x| < L - Bs \Rightarrow |\hat{f}(0) - f(x,y)| \leq C_b \delta/s$.\\

\noindent $\mathrm{(b)}~$ If $Q \subset A$ is measurable and invariant under rotations in the $y$-coordinate then the average value of $f$ on $Q$ is $\hat{f}(0)$.\\

\noindent $\mathrm{(c)}~\|f\|_2^2 \leq C\left[\|Ls \hat{f}(0)^2 + \mathrm{Im}(f)\|_2^2\right] $.\end{lemma}

\noindent Proof: The functions $\{z^n| n \in \zz\}$ are orthogonal, but not orthonormal, for the $L^2$ norm coming from any finite rotation-invariant measure; in particular they are orthogonal with respect to $q$-area on the cylinder. Therefore, by uniform convergence there is a Parseval identity on the region $|\log(|z|)| < L/s - \epsilon$ for holomorphic functions and the area measure $|q|$. Passing to the limit as $\epsilon \to 0$, one gets the Parseval identity:

$$\begin{array}{rl}\|\mathbf{1}_{A^*} \cdot f\|_2^2 & = \sum\limits_{n = -\infty}^\infty \int\limits_{-L}^L\int\limits_0^{2\pi s} |\hat{f}(n)|^2 |e^{n(x+iy)/s}|^2 \, dy \, dx \\ {} & = 2\pi s \left[ 2L|\hat{f}(0)|^2 + \sum\limits_{n = 1}^\infty \left(|\hat{f}(n)|^2 + |\hat{f}(-n)|^2\right) \frac{\sinh(2nL/s)}{n/s}\right].\end{array}$$

\noindent Unfortunately, the real/imaginary parts of $\hat{f}(n)z^n$ and $\hat{f}(m)z^m$ are only guaranteed to be orthogonal if $|m| \neq |n|$; however, if $L/s$ is sufficiently large, they are nearly orthogonal even if $m = -n$. More precisely, if $L > 0$, then for all $n \in \nn$ we have $$\frac{|\langle\mathrm{Im}(\hat{f}(n)z^n), \mathrm{Im}(\hat{f}(-n)z^{-n})\rangle|}{\|\mathrm{Im}(\hat{f}(n)z^n)\|_2\cdot \|\mathrm{Im}(\hat{f}(-n)z^{-n})\|_2} < c_{L/s} < 1,$$ \noindent and $c_{L/s} \to 0$ as $L/s \to \infty$, and the cost of pretending that they are in fact orthogonal is only a multiplicative constant which tends to 1 exponentially quickly as $L/s \to \infty$. (The intuitive explanation of this is that most of the norm of  concentrates at opposite ends of the cylinder $A^*$, so they essentially behave like functions with disjoint support; it is easy to estimate the contributions from regions where $|x| \geq 0$ and where $|x| \leq 0$.) This proves (c).\\

\noindent (b) follows from Fubini's theorem, since the series converges uniformly on each circle and only the term for $n = 0$ has non-zero average.\\

\noindent In all estimates that follow, the cost of pretending that these components are actually orthogonal is only a small multiplicative constant. Similarly, for sufficiently large $L$ we may replace $\sinh$ with the exponential function, at the cost of a multiplicative factor close to $1$. In particular, we have $$\|\mathbf{1}_{A^*}\mathrm{Im}(f)\|_2^2 \, \dot{\asymp} \, sL |\hat{f}(0)|^2 + s^2 \sum_{n = 1}^\infty \left(|(\hat{f}(n))|^2 + |\hat{f}(-n)|^2\right)\frac{e^{2 n L/s}}{n}.$$

\noindent Since the $L^2$ norm of $\mathrm{Im}(f)$ is bounded by $\delta^2$, this implies that for all $n > 0$, we get the following bound on Laurent coefficients: $\noindent |\hat{f}(\pm n)| \leq C \sqrt{n} (\delta/s) e^{-n L/s}.$

$$\begin{array}{rl}
    |f(x,y) - \hat{f}(0)| & \leq \sum\limits_{n = 1}^{\infty} C \sqrt{n} \frac{\delta}{s} e^{-n L/s} e^{n(L/s - B)}\\
    {} & \leq C \frac{\delta}{s} \sum\limits_{n = 1}^\infty \sqrt{n} e^{-nB}\\ {} & \leq C \delta/s.\\
\end{array}$$

\noindent The last inequality, and therefore part (a), are due to convergence of the infinite series. $\Box$\\

\section{Genus 2 example}

\noindent In this section we discuss a low-complexity example that contains the main analytic ideas, with simplified geometric input. The remaining sections will discuss a variant of Rafi's thick-thin decomposition that generalizes the decomposition of the surface into pieces that appear in this example. This decomposition is strongly motivated by the product regions theorem of Minsky \cite{ProdReg}, and the concept of active interval for projections to the factors in Minsky's theorem due to Rafi. The notion of active interval can be found in \cite{Fellow}, Theorem A. We postpone the careful discussion of these topics until we need them, but for those who already familiar with coarse geometry of Teichm\"uller space, we give a brief explanation of our choice of examples.\\

\noindent Let $(X,q)$ be formed as follows: let $T^*$ be a rectangular torus of side length $S$ with a geodesic slit of length $2\pi s$, and let $A^*$ be an $L \times s$ rectangle. Form a translation surface by identifying the sides of length $2L$ with each other and the sides of $A^*$ of length $s$ with sides of the slit. The two endpoints of the slit are identified, forming a cone point of angle $6\pi$. We will assume $s/S$ and $s/L$ are smaller than $b/2$, but allow them to be arbitrarily small.\\

\noindent With respect to the cylinder coordinates in the previous section, the average value of $f$ on the region $|x| < L - Bs$ is $\hat{f}(0)$. We now break up $X$ into regions on which we can apply familiar estimates. The regions are as follows:

\begin{itemize}\item $T$ will be the complement of the $bS$-neighborhood of the slit in the slit torus  $T^*$.
\item Let $U$ be the $(b \cdot s)$-neighborhood of the cone point.
\item $A$ will be the annulus $|x| < L - Bs$ contained in $A^*$
\item Let $E$ be the $Bs$-neighborhood of the slit with $U$ deleted.
\item Let $M$ be the remaining expanding annulus consisting of points distance between $Bs$ and $bS$ away from the slit.
\end{itemize}

\noindent Now, let $d = \max(s/S, \sqrt{s/L}).$  Let $\delta = \hat{\delta}d.$ We will show that if $\hat{\delta} << 1$, then it is impossible to have a meromorphic $f$ on $X$ with $\|f\|_2 = 1$ and $\|\mathrm{Im}(f)\|_2 < \delta.$\\

\noindent Pick a base point $p_T \in T.$ We run the following steps:\\

\noindent STEP 1: If $p_1,p_2 \in T$ then $|f(p_1) - f(p_2)| < \frac{C\delta}{S}.$\\

\noindent Proof: This follows from \Cref{Rodrigo} estimates, since the diameter and injectivity radius are uniformly bounded below by a multiple of $S$ in $T$, and the diameter of $T$ is bounded above by a multiple of $S$, and the value of $\|f\|_2$ is bounded by $1$.\\

\noindent STEP 2: $\hat{f}(0) < C/\sqrt{Ls}$.\\

\noindent Proof: By the Parseval identity $4\pi Ls|\hat{f}(0)|^2 \leq \|f\|_2^2 = 1$.\\

\noindent STEP 3: If $(x,y) \in A$ then $|f(x,y) - \hat{f}(0)| < C\delta/\sqrt{Ls}.$\\

\noindent Proof: This is exactly part (a) from \Cref{CylinderConstant}.\\

\noindent STEP 4: $|\mathbf{1}_T \cdot f|_\infty < C/S.$\\

\noindent Proof: If not, then $f >> 1/S$ on all of $T$, by step 1, and $T$ has area $S^2$, so $\int_T |f|^2 >> 1$. This contradicts $\|f\|_2 = 1.$\\


\noindent STEP 5: If $p_1, p_2$ are in the same expanding annular region (i.e. both in $M$ or both in $U$), then $|f(p_1)-f(p_2)| \leq \frac{C\delta}{\min [d_q(z_0,p_1),d_q(z_0,p_2)]}$.\\

\noindent Proof: Follows from integrating the \Cref{Rodrigo} estimate.\\

\noindent STEP 6: If $p_1,p_2 \in E, |f(p_1) - f(p_2)| < C\delta/s$.\\

\noindent Proof: This is essentially the same as step 1, except the geometry of $E$ is controlled by $s$ instead of $S$.\\

\noindent STEP 7: If $p_1, p_2 \in A \cup E \cup T \cup M$, then $|f(p_1) - f(p_2)| < C\delta/s$.\\

\noindent Proof: Follows easily from steps 1-6 and the triangle inequality.\\

\noindent STEP 8: If $\hat{\delta} << 1$, then we get a contradiction by showing that the sum of $f$ and a constant function has small $L^2$ norm, when by hypothesis the norm of such a function must be at least $1$.\\

\noindent For this step we have two cases, depending on whether the torus or cylinder has greater area (up to a factor of $4\pi$). In each case, we will take the function $f$ and subtract a constant function (coming from the larger piece) and show that the result has small $L^2$ norm. This contradicts the fact that $f$ has norm $1$ and is orthogonal to the constant functions.\\

\noindent CASE 1: $s/S \geq \sqrt{s/L}$, or equivalently $1 ~ \dot{\asymp} ~ Ls \geq S^2$. We thus have $\delta = \hat{\delta}s/S.$ Let $g(x) = f(x) - \hat{f}(0)$.\\

\noindent In this case we get the following:

\begin{center}

UPPER BOUNDS FOR $g = f - \hat{f}(0)$:\\

$~$\\

\begin{tabular}{ |c|c|c|c|c| } 
\hline 
 & & & & \\
$V$ & $Area(V) ~ \dot{\asymp} ~ *$ & $\left| ~ g\downharpoonright_{\partial V}\right| ~  \dot{\prec} ~ *$ & $\left|\int_V g |q| \right| \dot{\prec} ~  *$ & $\int_V |g|^2|q| ~\dot{\prec}~ *$\\
 & & & & \\
\hline
 & & & & \\
$A$ & $1$ & $\delta/s~\dot{\asymp}~\hat{\delta}$ & $\int_A g |q| =0 $ & $\hat{\delta}^2$ \\
 & & & & \\
\hline
 & & & & \\
$E$ & $s^2$ & $\delta/s = \hat{\delta}/S$ & $s\delta~\dot{\prec}~S\hat{\delta}$ & $\delta^2 ~\dot{\prec}~ \hat{\delta}^2$\\
 & & & & \\
\hline
 & & $\infty$ & & \\
 $U$ & $s^2$ & $~$ & $s\delta~\dot{\prec}~S\hat{\delta}$ & $\hat{\delta}^2$ \\
 & & $\delta/s~\dot{\prec}~\hat{\delta}/s$ & & \\
\hline
 & & $\delta/s = \hat{\delta}/S$ & & \\
 $M$ & $S^2$ & & $S\hat{\delta}$ & $\hat{\delta}^2$\\
 & & $\delta/S$ & & \\
\hline
 & & & & \\
$T$ & $S^2$ &$\delta/s = \hat{\delta}/S$  & $S\hat{\delta}$ & $\hat{\delta}^2$ \\
 & & & & \\
\hline
\end{tabular}
\end{center}

\noindent All terms in the second column are $O(\hat{\delta})$ but their sum is at least $1$, a contradiction.\\

\noindent CASE 2: $s/S \leq \sqrt{s/L}$, or equivalently $Ls \leq S^2 ~ \dot{\asymp} ~ 1$. Now, $\delta = \hat{\delta}\sqrt{s/L}$. Now, let $g(x) = f(x) - f(p_T)$. We make a similar upper bound table as before, but we suggestively switch the positions of $A$ and $T$.

\begin{center}

UPPER BOUNDS FOR $g = f - f(p_T)$

\begin{tabular}{ |c|c|c|c|c| } 
\hline 
 & & & & \\
$V$ & $Area(V) ~ \dot{\asymp} ~ *$ & $\left| ~ g \downharpoonright_{\partial V \setminus z_0}\right| ~  \dot{\prec} ~ *$ & $\left|\int_V g |q| \right| \dot{\prec} ~  *$ & $\int_V |g|^2|q| ~\dot{\prec}~ *$\\
 & & & & \\
\hline
 & & & & \\
$T$ & $1$ & $\delta~\dot{\asymp}~\hat{\delta}$ & $\delta~\dot{\asymp}~\hat{\delta}$ & $\delta^2~\dot{\prec}~\hat{\delta}^2$ \\
 & & & & \\
\hline
 & & & & \\
$E$ & $s^2$ & $\delta/s~\dot{\prec}~\hat{\delta}/s$ & $s\delta~\dot{\prec}~s\hat{\delta}$ & $\delta^2~\dot{\prec}~\hat{\delta}^2$\\
 & & & & \\
\hline
 & & $\infty$ & & \\
$U$ & $s^2$ &  & $s\hat{\delta}$ & $\hat{\delta}^2$ \\
 & & $\delta/s \dot{\prec} \hat{\delta}/s$ & & \\
\hline
 & & $\delta/s$ & & \\
 $M$ & $1$ & & $\delta$ & $\delta^2 ~\dot{\prec}~ \hat{\delta}^2$\\
 & & $\delta$ & & \\
\hline
 & & & & \\
$A$ & $Ls$ & $\delta/s = \hat{\delta}/\sqrt{Ls}$ & $\delta L = \hat{\delta} \sqrt{Ls}$ & $\delta^2  + \hat{\delta}^2~\dot{\asymp}~\hat{\delta}^2$ \\
 & & & & \\
\hline
\end{tabular}
\end{center}

\noindent This gives us a similar contradiction. Perhaps the final row of the table deserves some explanation. The entire goal is to control the value of the last entry using the Parseval identity from part (c) of \Cref{CylinderConstant}. We actually already know that the contribution to the $L^2$ norm of $f - \hat{f}(0)$ coming from $A$ is controlled by $\delta$, from the proof of part (c) of \Cref{CylinderConstant}. Note that the terms $\hat{g}(n)$ and $\hat{f}(n)$ agree for $n = 0$, so the contribution of this part is controlled by $C\delta^2$. Thus we only need to control $\hat{g}(0)$, which is done using the bounds boundary value, upgrading to a uniform bound via the maximum modulus principle, and multiplying this bound squared by area to estimate the integral.

\section{Thick-Thin Decompositions}

\subsection{Rafi's Thick-Thin Decomposition}

\noindent In this section we recall some results of Minsky and Rafi regarding the geometry of quadratic differentials. Throughout, a \emph{curve} will always mean a homotopy class of un-oriented simple closed curve in $X$ which is not homotopic to a puncture or a point.

\begin{definition}Let $b \in \rr^+, A_b = \{(x,y) \in \rr^2: 0 < y < b\}/(x,y) \thicksim (x + 1,y)$. If $A$ is an annulus conformal to $A_b$ we say that $b$ is the modulus of $A$, and write $\mathrm{Mod}(A) = b$.\end{definition}

\noindent Every annulus is conformal to either $\cc \setminus \{0\}$ or to $A_b$ for a unique $b$.

\begin{definition}If $\alpha$ is a curve in a Riemann surface $X$, and $A$ is an annulus in $X$, write $A \simeq \alpha$ if $A$ has a deformation retract in the homotopy class $\alpha$. We say $\alpha$ is the core curve of $A$. The \emph{extremal length} of $\alpha$ is defined to be $$\mathrm{Ext}_{\alpha}(X) := \left(\sup\limits_{A \simeq \alpha} \mathrm{Mod}(A)\right)^{-1}.$$ \end{definition}

\begin{definition}Let $(X,q)$ be a quadratic differential, and let $\alpha$ be a simple closed curve in $X$. We say $\alpha$ is a \emph{cylinder curve} if there is a $q$-geodesic representative of $\alpha$ that does not pass through a singular point of the $q$-metric.\end{definition}

\noindent If $\alpha$ is a cylinder curve, then there is a one-parameter family of geodesic representatives of $\alpha$ that do not pass through singular points, and the union of all such geodesics is an annulus in $X$, which we call a \emph{cylinder}. Moreover, these curves, plus the boundary curves of the cylinder are all the $q$-geodesic representatives of $\alpha$.\\

\noindent If $\alpha$ is not a cylinder curve, then in the completion $\overline{X}$ of $X$ (with respect to the $q$-metric), any length-infimizing sequence of constant speed parametrized loops in the class $\alpha$ converge (as parametrized geodesics, up to the action of the dihedral group) to a unique loop in $\overline{X}$, which we call the $q$-geodesic representative of $\alpha$.

\begin{definition} If $\alpha$ is a curve, we say that an (open) annulus $A \subset X$ is an \emph{expanding annulus for} $\alpha$ if the following hold:
\begin{itemize}
    \item $A$ has core curve $\alpha$
    \item $A$ If $\alpha$ is a cylinder curve, $A$ does not contain any geodesic in the cylinder of $\alpha$
    \item One boundary component of $A$ is a $q$-geodesic representative for $A$
    \item In the (completion of the) annular cover of $X$ corresponding to $A$, the two boundaries of $A$ are a uniform distance apart
    \item $A$ is maximal with respect to the previous properties.
\end{itemize}

\noindent As a convention, we will also say that a puncture $p$ admits an expanding annulus, which is the maximal singularity-free annulus whose boundary is equidistant from $p$ that deformation retracts to a loop about $p$. (The double cover of such an annulus will isometric to a finite cover of $\{z \in \cc: \}$ in its intrinsic geometry inherited from the $q$-metric)

\end{definition}

\noindent A curve can have at most two expanding annuli and at most one (maximal) cylinder. Since our surfaces are orientable, it makes sense to talk about which side of a curve an expanding annulus belongs to.  
We will use the following approximation to the modulus that can be computed more easily from flat geometry:

\begin{notation} If $A$ is a cylinder, let $\mu(A)$ be its modulus; if $A$ is an annulus not intersecting a cylinder whose boundaries are uniform distance apart, let $\mu(A)$ be the log of the ratio of the lengths of its boundary components (choose the ratio that is greater than 1 before taking the log). If $A$ only has one boundary curve and the other boundary is a puncture, set $\mu(A) = \infty$.\end{notation}

\noindent The following is a consequence of estimates from \cite{Harmonic}, section 4, reformulated in \cite{CRS}, section 5:

\begin{theorem} There is a number $\epsilon_0$, and positive constants depending only on $g$ and $n$, such that if $q$ is a quadratic differential on $X \in \mathcal{M}_{g,n}$, then if $\mathrm{Ext}_\alpha(X) < \epsilon_0,$ then $X$ contains an expanding annulus or cylinder $A$ with core $\alpha$ with $$\mu(A) ~ \dot{\asymp} ~ \frac{1}{\mathrm{Ext}_\alpha(X)}.$$ Moreover, of the expanding annuli and cylinder that exist for $\alpha$, $A$ is the one maximizing $\mu$. If $A$ is an expanding annulus, then there is a singularity-free annulus of $A^\prime \subset A$ whose boundaries are each constant distance from the boundaries of $A$, with $$\mu(A^\prime) ~ \dot{\asymp} ~ \mu(A).$$ $A^\prime$ is said to be a \emph{primitive regular annulus}.\end{theorem}

\noindent Now, fix a number $\epsilon_0$ smaller than the Margulis constant and small enough to meet the hypotheses of the theorem, and consider the collection of subsurfaces obtained by designating representatives of curves $\alpha$ of extremal length $\epsilon_0$ or smaller by deleting their geodesic representatives in the hyperbolic metric in the conformal class $X$.\\

\noindent Now, for each such subsurface $Y$, take the cover $\hat{Y}$ of $X$ associated to $\pi_1(Y)$, pull back the $q$-metric to this cover and call it the $\hat{q}$-metric. There is an open annulus of infinite area corresponding to each component of $\partial Y$, with geodesic boundary. If $\alpha$ is not a cylinder curve, there is a unique such annulus $A_\alpha$. If $\alpha$ is a cylinder curve, let $A_\alpha$ be the unique such annulus containing the entire cylinder. Delete from $\hat{Y}$ the union of the open annuli $A_\alpha$; we then say that $$\hat{Y} \setminus \bigcup\limits_{\alpha \subset \partial Y} A_\alpha =: Y_q$$ is the $q$-\emph{geodesic representative of the subsurface} $Y$.\\

\noindent Often, $Y_q$ is a surface with boundary, but there are degenerate cases in which it is not. In addition, some points on the boundary components of $Y_q$ may be identified in $X$. However, for any sufficiently small $\epsilon > 0$, the $\epsilon$-neighborhood of $Y_q$ in $\hat{Y}$ is homeomorphic to the corresponding hyperbolic surface with boundary $Y$. An example in section 5 of \cite{ThickThin} shows that $Y_q$ can be a spine of $Y$. However, we have the following definition and theorem from \cite{ThickThin} controlling the geometry of $Y_q$:

\begin{definition}Let $(X,q) \in QD(\mathcal{T}_{g,n})$ and let $Y_q$ be the $q$-geodesic representative of a thick component $Y$ of the thick-thin decomposition of $(X,q)$. If $Y$ is a pair of pants, let $\mathrm{size}_q(Y)$ be the maximum of the $q$-lengths of the boundary components of $Y$. Otherwise, let $\mathrm{size}_q(Y)$ be the minimum of the $q$-lengths of essential simple closed curves in $Y$. We say that $\mathrm{size}_q(Y)$ is the \emph{size} of the subsurface $Y$ in the $q$-metric.\end{definition}

\begin{theorem}\label{flatthickthin}Let $X,q,Y$ be as above. Let $\beta$ be any essential curve in $Y$, and let $\ell_q(\beta)$ denote its length in the $q$-metric and $\ell_\sigma(\beta)$ its length in the hyperbolic metric on $X$. Let $\diam_q$ and $Area_q$ denote diameter and area (respectively) in the $q$-metric. Then we have the following coarse estimates of the geometry of $Y_q$, in which all constants depend only on $g$ and $n$:

\begin{itemize}
    \item $\ell_q(\beta) ~ \dot{\asymp} ~ \mathrm{size}_q(Y) \ell_\sigma(\beta)$
    \item $\diam_q(Y_q) ~ \dot{\asymp} ~ \mathrm{size}_q(Y)$
    \item $Area_q(Y) := \int_{Y_q}|q| \leq C \cdot \mathrm{size}_q(Y)^2.\\$
\end{itemize}
\end{theorem}

\noindent We have the following lemma (\cite{EMR}, Lemma 3.5) relating size to expanding annulus:

\begin{lemma} There exist positive constants $C_1, C_2$, depending only on topology ($g$ and $n$), such that every thick component $Y$ of $(X,q)$ has the following property: every puncture or boundary component of $Y$ with $q$-geodesic length less than $C_1 \cdot \mathrm{size}_q(Y)$ admits an expanding annulus of circumference at least $C_2 \cdot \mathrm{size}_q(Y)$ in the direction of the subsurface $Y$.\end{lemma}

\begin{corollary}There is a constant $c > 0$, depending only on $g,n$, such that the following holds. If $Y$ is a thick component of $(X,q)$ and $Y$ is not a pair of pants, and every boundary curve $\alpha \subset \partial Y$ satisfies $\ell_q(\alpha) < c ~ \cdot ~ \mathrm{size}_q(Y)$, then $Area_q(Y) ~ \dot{\asymp} ~ \mathrm{size}_q(Y)^2.$ (If $Y$ is a pair of pants, then the same statement is vacuously true.) \end{corollary}

\noindent For non-pants, one direction is due to Rafi, and the other follows immediately if we can show any of the expanding annuli are contained in $Y$. Now, the $q$-diameter of the union of any two boundary curves is bounded below by a constant multiple of $\mathrm{size}_q(Y)$, since the union of two such curves and any arc connecting them contains a curve that is essential in $Y$. In particular, the $q$-geodesic representatives of the boundary components must be disjoint if $c$ is sufficiently small. This implies the expanding annuli in the direction of $Y$ do not leave $Y$. $\Box$

\subsection{The Primitive Annuli Decomposition}

\noindent We would like to use a modified thick-thin decomposition of a surface $(X,q)$, which consists of some annuli of large modulus and the complement of their union. We will not use maximal cylinders and expanding annuli exactly, however; we remove a bounded amount of modulus near each boundary component, so that each annulus satisfies the hypothesis of \Cref{shells} or \Cref{CylinderConstant}.\\

\begin{lemma}\label{SpecificAnnuli} There exist constants $\mu_0 > B > 2$ depending only on the topology of $X$, such that the following hold:\\

\noindent Suppose $(X,q)$ is taken from the principal stratum and then punctured at all singularities. For every expanding annulus $A$ of $(X,q)$ associated to a curve of puncture of resulting surface satisfying $\mu_0 < \mu(A) \leq \infty$, we pick two level circle of the expanding annulus: let $o_A$ be the (topological) circle $\gamma_A$ consisting of points $p$ such that the distance from $p$ to the the puncture or geodesic boundary component of $A$ is $1/B$ times the distance from $p$ to the the other boundary of $A$ in the completion $\overline{X}$ of $X$ (with respect to the $q$-metric). If $A$ is homotopic to a puncture let $i_A$ be the puncture, and otherwise let $i_A$ be the simple closed curve contained in $A$ that consists of points $B$ times the $q$-length of the $q$-geodesic boundary component of $A$. For each cylinder $A$ with $\mu(A) > \mu_0$, let $i_A$ and $o_A$ be geodesics homotopic to the core curve of $A$, each distance $B$ times the $q$-length of the core curve away from a boundary component.\\

\noindent Then, the curves $i_A$ and $o_A$ are pairwise disjoint for each annulus $A$, and $i_A$ is shorter unless they are both core curves of a cylinder. Moreover, the corresponding annuli are disjoint from each other; that is, for $A_1, A_2$ distinct, the annulus bounded by $i_{A_1}$ and $o_{A_1}$ is disjoint from the annulus bounded by and $i_{A_2}$ and $o_{A_2}$.\\

\noindent The constant $B$ can be taken arbitrarily large.
\end{lemma}

\noindent Proof: To make use of non-positively curved geometry we pass to the completion of the orienting double cover of $(X,q)$. The universal cover of this completion will be denoted $(\hat{X},\hat{q})$. It is complete, homeomorphic to a disk, and non-positively curved in the sense of Alexandrov.\\

\noindent If $S_q$ is the set of singularities of $q$, note that in the orienting double cover, the inverse images of any expanding annulus $A$ corresponding to an element of $S_q$ or a curve, must satisfy one of the following three possibilities:\\

\noindent 1. It consists of a single annulus which is a double cover of $A$, and it is the expanding annulus associated to the inverse image of the corresponding puncture or curve, which is a single curve.\\

\noindent 2. It consists of a pair of disjoint annuli, which are the expanding annuli for a pair of disjoint curves in the orienting double cover.\\

\noindent 3. It consists of a pair of disjoint annuli $A_1$ and $A_2$, each of which consists of a disjoint union of level circles for distance to a geodesic lift of the geodesic boundary of $A$, but which are proper subsets of disjoint expanding annuli for two curves that are the preimage of the core curve of $A$, and the outer boundaries of $A_1$ and $A_2$ intersect some boundary $A_2$ and $A_1$, respectively. (This is necessary for $A$ to be a maximal expanding annulus.)\\

\noindent In case 3, we claim that the width of of $A_i$ is more than a third of the width of the canonical expanding annulus containing $A_i$, as long as $\mu_0$ is large enough to guarantee that the width of $A$ exceeds the length of the geodesic boundary of $A$. (By theorem 4.5 of \cite{Harmonic} such a $\mu_0$ exists.) Indeed, such an annulus about $A_1$ would intersect both boundary components of $A_2$, and then expand an additional distance of at least the length of the lifts of the $q$-geodesic of $A$, which would imply that it contained a singularity.\\

\noindent Clearly, it is enough for us to prove that the preimages of these annuli to the orienting double cover are all disjoint, or more generally, that the connected components of the preimages in $(\hat{X},\hat{q})$ are all disjoint.\\

\noindent The annulus bounded by $i_A$ and $o_A$ is covered by a disjoint union of connected components, which are covers of degree $1,2$, or $\infty$ (with deck group isomorphic to $\zz$) in $(\hat{X},\hat{q}).$ Now, any such connected component consists of a pair $i_A$ and $o_A$. Any connected component $\gamma$ of the inverse image of a level circle of $A$ in $(\hat{X},\hat{q})$, when deleted, divides $(\hat{X},\hat{q})$ into two components, since such a component is either a topological circle, or fellow travels a bi-infinite geodesic. (For a complete, Alexandrov non-positively-curved metric on the disk, there is a boundary at $\infty$ homoemorphic to a circle, and we can just collapse this boundary and invoke the Jordan separation theorem.) Denote by $I_{\gamma}$ and $O_\gamma$ the two parts of the induced partition of the singularities of $(\hat{X},\hat{q})$, by which we mean the collections of inverse images of singularities in $(X,q)$ lying in the two connected components, \emph{including those which have cone angle} $2\pi$. Distinguish between the two sets in the following way: the singularities on $I_{\gamma}$ are those that are closer the corresponding connected component of the preimage of $\hat{i}_A$ of $i_A$ than to the component $\hat{o}_A$ that maps to $o_A$, and the singularities in $O_\gamma$ are closer to $\hat{o}_A$ than to $\hat{i}_A$.\\

\noindent First we deal with the case in which one of the annuli is a cylinder. Clearly, if two cylinders intersect, then their core curves intersect essentially, which contradicts the fact that they both admit cylinders of modulus greater than $\mu_0$. However, the width of an primitive expanding annulus that contains a point on the boundary of a cylinder is at most the length of the core curve of the cylinder, because it does not contain two distinct preimages of that point. Therefore it does not contain any point in the cylinder which is distance $B$ times the length of the core curve.\\

\noindent Now, we observe that $I_{A_1} = I_{A_2}$ if and only if they arise from the same lift of an annulus, since their core curves are homotopic if finite, and they lie on the same side(s) of the same infinite geodesic or cylinder of geodesics if $A_1$ and $A_2$ are infinite.\\

\noindent Now, we break the analysis of pairs of expanding annuli into three cases.\\

\noindent CASE 1: $I_{A_1} \cap I_{A_2} = \emptyset$. If $B$ is large enough, every point in the annulus or strip bounded by $\hat{i}_{A_1}$ and $\hat{o}_{A_1}$ is closer to $I_{A_1}$ than to $I_{A_2}$. We have a similar statement if we reverse the indices $1$ and $2$, and combining these gives the desired disjointness.\\

\noindent CASE 2: $I_{A_1}$ and $I_{A_2}$ are neither disjoint nor nested. Let $\ell_i$ be the lengths of the geodesic corresponding to $A_i$ and let $W_i$ be the widths. Then the curve with length $\ell_1$ contains points in $I_{A_2}$ and $O_{A_2}$, so $\ell_1 > W_2$. Similarly $\ell_2 > W_1$. We thus get $\ell_1 + \ell_2 > W_1 + W_2$. But for all sufficiently large $\mu$ we have $\ell_i < W_i$, so this is impossible.\\

\noindent CASE 3: They are nested; without loss of generality $I_{A_1} \subsetneq I_{A_2}$. Let $\ell_i, W_i$ be as in the previous case. Then $A_2$ cannot correspond to a single point, and is therefore covered by components of infinite diameter. Consider the bi-infinite geodesic corresponding to $A_2$. Now, a component of the inverse image of $A_1$ begins at a singularity in $I_{A_2}$ or a geodesic joining a bi-infinite sequence of singularities inside of $I_{A_2}$, and it can expand in the direction containing $O_{A_2}$ for an additional width of at most $\ell_2$ after it first meets the geodesic boundary of the lift of $A_2$, just as in our analysis when one of $A_1$ and $A_2$ was a cylinder. Just as in that case, we conclude that the expanding annulus for $A_1$ does not reach far enough past the geodesic boundary if $B > 1$. $\Box$

\begin{lemma}\label{CorrectSize} Let $\Sigma$ be the set of singularities of $(X,q)$. If $\mu_0$ is chosen sufficiently large in \Cref{SpecificAnnuli}, the cylinders and expanding annuli from \Cref{SpecificAnnuli} are deleted from $(X,q)$, and $\Sigma$ is also deleted, then for each remaining connected component $Y_i$, there is a corresponding component $Z_i$ of the thick-thin decomposition of $(X \setminus \Sigma, q)$ representing the homotopy class of the subsurface $Y_i$. Let $p \in Y_i$, and $r(p)$ denote the radius of the largest embedded disk in $X \setminus \Sigma$ centered at $p$. Let $\gamma$ be any boundary component of $Y_i$, with length $\ell(\gamma)$.\\

\noindent We then have the following estimates: $\ell(\gamma),r(p),\sqrt{Area_q(Y_i)} ~ \dot{\asymp} ~ \mathrm{size}_q(Z_i)$.\\

\noindent All implied constants depend only on $g, n, B, ~\mathrm{and}~ \mu_0.$\end{lemma}

\noindent We note that it does not matter whether we measure these with respect to the intrinsic metric on $Y_i$ or on the ambient metric on $X$.\\

\noindent Proof: This follows easily from Rafi's thick-thin decomposition and \Cref{SpecificAnnuli}.$\Box$

\begin{definition}Let $\mathcal{S}$ vary over all systems $U$ of curves from \Cref{SpecificAnnuli} that separate $(X,q)$ into two disjoint components, neither of which is a disjoint union of expanding annuli. If $X_U$ denotes the component with smaller area, and $\ell(U)$ denotes the sum of the lengths of the curves in $U$, define the \emph{hourglass ratio} of $(X,q)$ by $$H(X,q) := \min \left(\{1\} \cup \left\{ \frac{\ell(U)}{\mathrm{Area}_q(U)^{1/2}}: U \in \mathcal{PAD} \right\}\right).$$\end{definition}

\noindent Remark: We call this the hourglass ratio because it measures to what extent there is a small passage separating two much larger components. We also remark that we do not require the two components to be connected, but if we did this would only change the value of $H(X,q)$ up to a multiplicative constant.\\

\begin{definition} We refer to the collection of components bounded by the collection of simple closed curves in \cref{SpecificAnnuli} as the \emph{Primitive Annuli Decomposition} of $(X,q)$ and denote is by $\mathcal(X,q)$\end{definition}

\section{Meromorphic Functions and Efficient Paths}

\begin{proposition}\label{EfficientPaths}Let $Y$ be a component of $\mathcal{PAD}(X,q)$, i.e. a connected component of the the complement of the specific collection of simple closed curves on $(X,q)$ from \Cref{SpecificAnnuli} and \Cref{CorrectSize}. Let $f$ be meromorphic and $L^2$ on $(X,|q|).$ If $Y$ is not an expanding annulus, then for any $p_1,p_2 \in Y$ we have $$|f(p_1) - f(p_2)| ~ \dot{\prec} ~ \|\mathrm{Im}(f)\|_2/\|\mathbf{1}_Y\|_2.$$ Moreover, if there is a path from $x_1 \in X$ to $x_2 \in X$ in $X$ which is distance at least $r$ away from every singularity, then $$|f(x_1) - f(x_2)| ~ \dot{\prec} ~ \|\mathrm{Im}(f)\|_2/r.$$ \end{proposition}

\noindent Proof: The first claim follows from \Cref{Rodrigo} and \Cref{CorrectSize} on all components except for cylinders. For cylinder components we can use \Cref{CylinderConstant}.\\

\noindent The second claim is not much different. We can take our original path and consider any components it enters. Any two points in the same thick component $Y_i$ are joined by a path in $Y_i$ whose length is $O(\mathrm{size}_q(Y_i))$. Any two points on an expanding annulus can be joined by an arc which consists of an arc constant distance from both boundary components and an arc that is perpendicular to all such arcs. We can assume that our path takes this form whenever it enters a non-cylinder component; since the distance to singularities is constant up to a bounded multiple on each such $Y_i$ this does not cause the closest approach to a singularity to decrease by more than a bounded factor.\\

Now, our new path enters each thick component and cylinder at most once. We integrate the estimate from \Cref{Rodrigo} on all excursions into expanding annulus components, and apply the first claim for segments of paths that belong to components that are not expanding annuli. $\Box$

\begin{proposition}Let $\|\cdot\|_2$ denote the $L^2$ norm with respect to the $q$-area on a half-translation surface $(X,q)$. There is a constant $C_{g,n} > 0$ depending only on the genus and number of marked points of $X$, such that if $f$ is a nonzero $L^2$ meromorphic function on a unit area half-translation surface $(X,q)$ with $\int\limits_X f |q| = 0,$ then $$\|\mathrm{Im}(f)\|_2 > C_{g,n}H(X,q)\|f\|_2.$$ \end{proposition}

\noindent Proof: To simplify we may assume $(X,q)$ has unit area and $\|f\|_2 = 1$ and $\|\mathrm{Im}(f)\|_2 = \hat{\delta}H(X,q)$. We get a contradiction if $\hat{\delta} << 1.$\\

\noindent As usual, we will use the letter $C,C^\prime,C^{\prime\prime}$ to denote various positive constants depending only on $g$ and $n$. Their values may change from step to step.\\

\noindent STEP 1: Some component $X_0$ of $\mathcal{PAD}(X,q)$ has area at least $C$, and we may take $X_0$ to not be an expanding annulus.\\

\noindent Proof: The number of components of $\mathcal{PAD}(X,q)$ is bounded depending only on $g$ and $n$, and the area of an expanding annulus component is at most $C$ times the area of its neighboring component(s). So we may assume that the largest cylinder or largest non-annular component has area at least $C > 0$; we may take this component to be $X_0$.\\

\noindent STEP 2: If $x_0 \in X_0$, and $X_1$ is any component of $\mathcal{PAD}$ which is not an expanding annulus, then $x_1 \in X_1$ then $\left|f(x_1) - f(x_0)\right| < C\hat{\delta} \mathrm{Area}_q(X_1)$.\\

\noindent Proof: Let $g(x) = f(x) - f(x_0).$ We will first prove this for some $x_1 \in X_1$, and then extend to all $x_1 \in X_1$. On each piece $X_1$ that is not an expanding annulus component, $f$ is constant up to an additive error of $C \frac{C\hat{\delta}}{\mathrm{size}_q(X_i)}$ by \cref{EfficientPaths}. Moreover, if $X_1$ is a component of $\mathcal{PAD}(X,q)$ but not an expanding annulus, then there is a path $\gamma$ from $X_0$ to $X_1$ such that $$\frac{d(\gamma, \Sigma)}{\sqrt{Area(X_1)}} > CH(X,q).$$ \noindent So for some $x_1 \in X_1,$ an application of \cref{EfficientPaths} gives us $$H(X,q)|g(x_1)| \dot{\prec} \frac{\|\mathrm{Im}(f)\|_2}{\mathrm{Area}_q(X_1)^{1/2}}.$$

\noindent Moreover, by \Cref{EfficientPaths}, this is actually true for \emph{all} $x_1 \in X_1$. The claim then follows by dividing through by $H(X,q)\mathrm{Area}_q(X_1)^{1/2}.$\\

\noindent STEP 3: Proof of the proposition. We now assume $\hat{\delta} << 1.$ As in step 2, set $g(x) = f(x) - f(x_0)$. Then we must have $\|g\|_2 \geq 1.$ since $\|f\|$ is orthogonal to constants and $f - g$ is constant. However, by Step 2 we clearly have $\int_{X_1} |g|^2|q| \dot{\prec} (\hat{\delta})^2.$\\

\noindent Step 2 also implies that $$\|g(x)\| \dot{\prec} \frac{\hat{\delta}}{d(x,\Sigma)}$$ on all components that are not expanding annuli. We may use this as the boundary condition needed to apply \Cref{shells}. Summing these we conclude $\|g\|_2 < C\hat{\delta}$, a contradiction. $\Box$\\

\noindent Proof of \cref{HodgeGap}: Follows immediately as as corollary given the discussion in section 2.$\Box$

\section{Contraction along axes of pseudo-Anosov homeomorphisms}

\noindent A pseudo-Anosov homeomorphism induces a map on the space of measured foliations, with north-south dynamics on the space $\mathcal{PMF}$ of projective measured foliations. If the axis of the pseudo-Anosov diffeomorphism is contained in the principal stratum of quadratic differentials, we may use the Hodge norm as a norm on the tangent space to the space of quadratic differentials, which is locally $H_{odd}^1(\tilde{X}; \cc).$ If a class is of the form $\alpha + i \beta$ with $\alpha, \beta \in H_{odd}^1(\tilde{X}$, the flow in period coordinates is given by $g_t(\alpha + i \beta) = e^{t}\alpha + i e^{-t}\beta.$ Let $\phi$ be pseudo-Anosov with translation length $T$. One can apply the flow and then the inverse of the pseudo-Anosov homeomorphism to get a self-map of $H_{odd}^1(X;\rr)$. As is well known, this map is symplectic, and the eigenvectors have norms in $[e^{-T}, e^T].$ The top and bottom eigenvalues are simple and positive real; they correspond to the classes of the horizontal and vertical foliations. This map is also symplectic; so the eigenvalues take the form $$e^T = \lambda_1 > \lambda_2 \geq ... \geq \lambda_2^{}-1 > \lambda_1^{-1} = e^{-T}.$$

\begin{theorem}Let $\phi$ be pseudo-Anosov of translation distance $T$. Let $(X,q)$ be a half-translation surface on the axis of $\phi$ with $(X_t,q_t) = g_t(X,q)$. If $\phi$ belongs to the principal stratum, then $$\log(\frac{\lambda_1}{\lambda_2}) \geq C_{g,n}\int_0^T H(X_t,q_t) dt. \Box$$ \end{theorem}

\noindent We would like to have a similar theorem when $\phi$ does not have axis in the principal stratum. In this case the period coordinates do not make sense, but one may find a neighborhood of $X$ covered by a finite collection of cones in vector spaces of the form $H_{odd}^1(X^\prime,q^\prime) \times H_{odd}^1(X^\prime,q^\prime),$ where $(X^\prime,q^\prime)$ is a nearby quadratic differential in the principal stratum. Then, possibly after passing to a power of $\phi$, we may assume that our pseudo-Anosov does not non-trivially permute these cones. Then we have self-maps of these cones, and the eigenvalues of these maps satisfy $$e^T = \lambda_1 > \lambda_2 \geq ... \geq \lambda_2^{}-1 > \lambda_1^{-1} = e^{-T}.$$

\noindent We give a brief sketch of the proof here. Relevant definitions can be found in \cite{Hypertrain}.\\

\noindent Step 1: We consider the coding of a geodesic in the principal stratum that fellow travels the axis of $\phi$, via a train track splitting sequence $T$, following \cite{Hypertrain}. Such an axis will contain a subsequence that is a splitting sequence for a power of $\phi$. Howeotver, n every edge splits, so we fix some periodic splitting sequence $S$ in which every edge splits, which can be concatenated with $T$. We let $L_T$ and $L_S$ be the associated maps on period coordinates. We consider the pseudo-Anosov homeomorphisms $\phi_k$ associated to the splitting sequences $T^k S$ as $k \to \infty$.\\

\noindent We make the following claims:\\

\noindent Claim 1: The axes of the $\phi_n$ all belong to a fixed compact subsect of $\mathcal{M}_{g,n}$ that depends only on $T$ and $S$.\\

\noindent Claim 2: For each $\epsilon > 0$, there is some $k_0$ such that for all $k > k_0$, all but a fraction $(1 - \epsilon)$ of the length of the axis of $\phi_k$ $\epsilon$-fellow travels the axis of $\phi_k$  in a parametrized fashion. (This follows easily from hyperbolic properties of the flow on compact invariant sets, see e.g. \cite{HCIS}).\\

\noindent Claim 3: The average value of the hourglass ratio along the axis of $\phi_k$ converges to the average value of the hourglass ratio along the axis of $\phi$.\\

\noindent Claim 4: If $\lambda_2$ and $\lambda_1$ are the top two eigenvalues for $L_U$, for $U = T,S, T^NS$ etc. Then the spectral gap of of a linear map on the cone $X$ of equivalence classes tangential measures, is $\frac{\lambda_2}{\lambda_1} = \lim\limits_{k \to \infty} \diam_{Hilb}(U^kS X)$, where $X$ is the cone of equivalence classes of tangential measures, and $\diam_{Hilb}$ is the diameter in the Hilbert metric on $X$. From \cite{Hypertrain} we know that some power of $L_S$ is a contraction on $X$ in the Hilbert metric, and if we linearize the Hilbert metric near the attracting fixed point on a hyperplane representing the cone $X$ up to scaling, then the derivative of the map $L_S$ on $X$ with respect to the Hilbert metric is $L_S/\lambda_1 \circ P$, where $P$ is a projection onto the hyperplane.\\

\noindent We then get the conclusion by taking square the $N^{\mathrm{th}}$ root and the limit as $N \to \infty$.\\

\noindent A construction of Bell and Schleimer in \cite{SlowNS} gives examples of pseudo-Anosov homeomorphisms with $\lambda_2/\lambda_1$ arbitrarily close to 1, whenever the complex dimension of the moduli space $3g-3+n$ is at least $4$. We remark that our construction implies that such homeomorphisms must live deeper and deeper in the cusps of moduli space or have quasiconformal dilatation tending to $1$; however, it is well known that the quasiconformal dilatation is bounded below by a constant greater than 1 for each fixed $g,n$.\\

\section{Invariant Transverse Measures and Unique Ergodicity}

\noindent The following theorem is from \cite{SmithUE}, building on a theorem from \cite{Rodrigo}, which proved the result in the case when the $(X,q)$ has $q = \omega^2$ for some holomorphic 1-form $\omega$.\\

\begin{theorem}Let $\kappa(t)$ denote the systole of $(X_t,q_t) = g_t(X,q)$ with respect to the $q_t$-metric. If $$\int_0^\infty \kappa(t)^2 dt$$ diverges, then the vertical foliation on $(X,q)$ is uniquely ergodic.\end{theorem}

\noindent We would like to replace $\kappa(t)^2$ with the square of the hourglass ratio, but such a theorem cannot be true. A counterexample can be given by Strebel differential. That is, a surface obtained from a rectangle by gluing the top and bottom by a translation, and the left and right sides by a piecewise-translation. However, in this case, the invariant ergodic measures are all topologically equivalent.\\

\noindent Following the methods of Smith and Trevi\~{n}o, one may prove

\begin{theorem}Let $(X_t,q_t) = g_t(X,q)$. If $$\int_0^\infty H(X_t,q_t)^2 \, dt$$ diverges, then the invariant transverse measures on vertical foliation represent only one point in $\mathcal{MF}$.\end{theorem}

\noindent Proof: It is well known that $(X,q)$ decomposes into cylinders of closed leaves of the vertical foliation and minimal components for the foliation, see for example the survey of Masur and Tabachnikov \cite{Masur2002RationalBA}.\\

\noindent If the integral diverges, we observe that the surface is either a Strebel differential (in which case the result is trivial) or the vertical foliation is minimal. Indeed, if not there are at least two cylinders or minimal components whose area does not tend to zero but which are separated by a system of curves whose length decays exponentially, so $H(X_t,q_t) = O(e^{-t})$.\\

\noindent In the minimal case, we may further reduce to the case in which $(X,q)$ is in the minimal stratum by a small deformation along the strongly stable leaf into the principal stratum to some $(X^\prime,q^\prime)$; then $H(X_t,q_t) - H(X_t^\prime,q_t^\prime)$ decays exponentially.\\

\noindent The argument of Trevi\~{n}o and Smith now carries through with the following modification: our estimate for the $L^2$ norm of the imaginary part of a meromorphic function orthogonal to the constant functions may be substituted into their argument.\\

\section{Further Questions}

\noindent We expect that there should be similar geometric criteria that give bounds for the derivative of the Hodge norm on bigger larger subspaces.  For example, we expect that mutually orthogonal $L^2$ meromorphic functions must essentially live on disjoint subsurfaces and cylinders, i.e disjoint Minsky product region factors.\\

\noindent More precisely, for each $n$, one may define a ``level $n$ hourglass constant" and use this to give upper and lower bounds the derivative of the Hodge norm on the $n^{\tH}$ exterior power of $H_{odd}^1(\tilde{X},\rr).$ This can be done by similar methods.

\noindent If such an integral diverges, we would hope to be able to conclude that the veritcal foliation of $(X,q)$ has ergodic measures forming a simplex with at most $(n-1)$ vertices in $\mathcal{PMF}$. This could be seen as a quantitative generalization of McMullen's generalization of Masur's criterion (\cite{StableMasur}, Thm. 1.4).\\

\noindent Similarly, one should also be able to bound more Lyapunov exponents of pseudo-Anosov maps by similar integrals.\\

\noindent Finally, we leave as a question the qualitative sharpness of our spectral gap estimate: that is, can it be improved in any way except by estimating the multiplicative constants?\\



\bibliographystyle{amsalpha}
\bibliography{HodgeTori}

\providecommand{\bysame}{\leavevmode\hbox to3em{\hrulefill}\thinspace}
\providecommand{\MR}{\relax\ifhmode\unskip\space\fi MR }
\providecommand{\MRhref}[2]{%
  \href{http://www.ams.org/mathscinet-getitem?mr=#1}{#2}
}
\providecommand{\href}[2]{#2}
\begin{thebibliography}{ABEM12}

\bibitem[ABEM12]{abem}
Jayadev Athreya, Alexander Bufetov, Alex Eskin, and Maryam Mirzakhani,
  \emph{Lattice {P}oint {A}symptotics and {V}olume {G}rowth on
  {T}eichm{\"u}ller spaces}, Duke Mathematical Journal \textbf{161} (2012),
  no.~6, 1055--1111.

\bibitem[BS15]{SlowNS}
Mark Bell and Saul Schleimer, \emph{Slow north-south dynamics on
  $\mathcal{PML}$}, Groups, Geometry, and Dynamics \textbf{11} (2015).

\bibitem[CRS08]{CRS}
Young-Eun Choi, Kasra Rafi, and Caroline Series, \emph{Lines of minima and
  {T}eichm{\"u}ller geodesics}, Geom. Funct. Anal. \textbf{18} (2008), no.~3,
  698--754.

\bibitem[EM11]{EM}
Alex Eskin and Maryam Mirzakhani, \emph{Counting {C}losed {G}eodesics in
  {M}oduli {S}pace}, J. Mod. Dyn. \textbf{5} (2011), no.~1, 71--105.

\bibitem[EMR19]{EMR}
Alex Eskin, Maryam Mirzakhani, and Kasra Rafi, \emph{Counting closed geodesics
  in strata}, Inventiones mathematicae \textbf{215} (2019), no.~2, 535--607.

\bibitem[FMZ12]{FMZ}
Giovanni Forni, Carlos Matheus, and Anton Zorich, \emph{Lyapunov spectrum of
  invariant subbundles of the {H}odge bundle}, Erg. Thy. Dyn. Sys. \textbf{34}
  (2012), no.~2, 353--408.

\bibitem[For01]{FE}
Giovanni Forni, \emph{Deviation of ergodic averages for area-preserving flows
  on surfaces of higher genus}, Annals of Mathematics \textbf{154} (2001),
  no.~1, 1--103.

\bibitem[Fra]{Hypertrain}
Ian Frankel, \emph{{C}{A}{T}(-1)-type properties for {T}eichm\"uller space},
  arXiv: 1808.10022 [math.GT].

\bibitem[Ham10]{HCIS}
Ursula Hamenst{\"a}dt, \emph{Dynamics of the {T}eichm{\"u}ller {F}low on
  {C}ompact {I}nvariant {S}ets}, Inventiones Mathematicae \textbf{4} (2010),
  no.~2, 393--418.

\bibitem[Ham13]{HB}
\bysame, \emph{Bowen's {C}onstruction for the {T}eichm{\"u}ller {F}low}, J.
  Mod. Dyn. \textbf{7} (2013), no.~4, 489--526.

\bibitem[Hub06]{Hubbard}
John Hubbard, \emph{Teichm{\"u}ller {T}heory and {A}pplications to {G}eometry,
  {T}opology, and {D}ynamics. {V}olume {I}: {T}eichm{\"u}ller {T}heory}, Matrix
  Editions, Ithaca, New York, 2006.

\bibitem[Mar70]{Margulis}
Grigory~Aleksandrovich Margulis, \emph{On {S}ome {A}spects of the {T}heory of
  {A}nosov {S}ystems}, Ph.D. thesis, 1970, Springer, 2003.

\bibitem[McM]{StableMasur}
Curtis~T. McMullen, \emph{Diophantine and ergodic foliations on surfaces},
  Journal of Topology \textbf{6}, no.~2, 349--360.

\bibitem[Min92]{Harmonic}
Yair Minsky, \emph{Harmonic maps, length, and energy in {T}eichm{\"u}ller
  space}, J. Diff. Geom. \textbf{35} (1992), 151--217.

\bibitem[Min96]{ProdReg}
\bysame, \emph{Extremal length estimates and product regions in
  {T}eichm{\"u}ller space}, Duke Mathematical Journal \textbf{83} (1996),
  no.~2, 249--286.

\bibitem[MT02]{Masur2002RationalBA}
Howard~A. Masur and Serge Tabachnikov, \emph{Rational billiards and flat
  structures}, Handbook of Dynamical Systems \textbf{1A} (2002), 1015--1089.

\bibitem[Raf07]{ThickThin}
Kasra Rafi, \emph{Thick-{T}hin {D}ecomposition for {Q}uadratic
  {D}ifferentials}, Math. Res. Let. \textbf{14} (2007), no.~2, 333--341.

\bibitem[Raf14]{Fellow}
\bysame, \emph{Hyperbolicity in {T}eichm{\"u}ller space}, Geometry and Topology
  \textbf{18} (2014), 3025--3053.

\bibitem[Smi]{SmithUE}
M.~E. Smith, \emph{On the unique ergodicity of quadratic differentials and the
  orientation double cover}, arXiv: 1704.06303 [math.DS].

\bibitem[Tre14]{Rodrigo}
Rodrigo Trevi{\~n}o, \emph{On the ergodicity of flat surfaces of finite area},
  Geom. Funct. Anal. \textbf{24} (2014), no.~1, 360--386.

\end{thebibliography}

\end{document}